\documentclass{article}
\usepackage[english]{babel}

\usepackage[letterpaper,top=2cm,bottom=2cm,left=2cm,right=2cm,marginparwidth=1.75cm]{geometry}

\usepackage{calrsfs}
\usepackage{amssymb}
\usepackage{amsmath}
\usepackage{amsthm}
\usepackage{subcaption}
\usepackage{multirow}
\usepackage{graphicx}
\usepackage{amsfonts}
\usepackage{mathtools}
\usepackage{mathrsfs}
\usepackage{placeins}
\usepackage{relsize}

\usepackage{booktabs}
\newtheorem{remark}{Remark}
\newtheorem{definition}{Definition}

\newtheorem{problem}{Problem}

\usepackage[mathscr]{euscript}
\usepackage[dvipsnames]{xcolor}
\usepackage[colorlinks=true, allcolors=OliveGreen]{hyperref}
\usepackage{cleveref}
\usepackage{filecontents}
\usepackage{tikz}
\usepackage{pgfplots}
\usetikzlibrary{external}
\usepackage{authblk}
\usepackage{todonotes}
\usepackage{array}
\newcolumntype{C}{>{\centering\arraybackslash}p{2.5cm}}

\usepackage{caption}
\usepackage{subcaption}
\usepackage{tikz}
\usepackage{csvsimple}
\usepackage{siunitx}
\usepackage{tikz} 
\usepackage{pgfplots}
\pgfplotsset{compat=newest} 
\usepgfplotslibrary{units} 
\newcommand{\logLogSlopeTriangle}[5]
{

    \pgfplotsextra
    {
        \pgfkeysgetvalue{/pgfplots/xmin}{\xmin}
        \pgfkeysgetvalue{/pgfplots/xmax}{\xmax}
        \pgfkeysgetvalue{/pgfplots/ymin}{\ymin}
        \pgfkeysgetvalue{/pgfplots/ymax}{\ymax}

        \pgfmathsetmacro{\xArel}{#1}
        \pgfmathsetmacro{\yArel}{#3}
        \pgfmathsetmacro{\xBrel}{#1-#2}
        \pgfmathsetmacro{\yBrel}{\yArel}
        \pgfmathsetmacro{\xCrel}{\xArel}

        \pgfmathsetmacro{\lnxB}{\xmin*(1-(#1-#2))+\xmax*(#1-#2)} 
        \pgfmathsetmacro{\lnxA}{\xmin*(1-#1)+\xmax*#1} 
        \pgfmathsetmacro{\lnyA}{\ymin*(1-#3)+\ymax*#3} 
        \pgfmathsetmacro{\lnyC}{\lnyA+#4*(\lnxA-\lnxB)}
        \pgfmathsetmacro{\yCrel}{\lnyC-\ymin)/(\ymax-\ymin)} 

        \coordinate (A) at (rel axis cs:\xArel,\yArel);
        \coordinate (B) at (rel axis cs:\xBrel,\yBrel);
        \coordinate (C) at (rel axis cs:\xCrel,\yCrel);

        \draw[#5]   (A)-- node[pos=0.5,anchor=north] {}
                    (B)-- 
                    (C)-- node[pos=0.5,anchor=west] {#4}
                    cycle;
    }
}

\DeclareMathAlphabet{\mathcalligra}{T1}{calligra}{m}{n}

\graphicspath{{Images/}}

\tikzset{  font={\fontsize{15pt}{12}\selectfont}}

\title{High-order Discontinuous Galerkin Methods for the Monodomain and Bidomain Models}

\author[1]{Federica Botta\footnote{These authors equally contributed to this work}\textsuperscript{,}}
\affil[1]{federica.botta@mail.polimi.it}

\author[2]{Matteo Calafà\textsuperscript{$\ast$,}}
\affil[2]{Institut for Mekanik \& Produktion, Aarhus University, Inge Lehmanns Gade 10, 8000 Aarhus, Denmark}

\author[3]{Pasquale C. Africa}
\affil[3]{mathLab, SISSA,
via Bonomea, 265,
34136 Trieste, Italy}

\author[4]{Christian Vergara}
\affil[4]{LABS, Dipartimento di Chimica, Materiali e Ingegneria Chimica, Politecnico di Milano, Piazza Leonardo da Vinci 32,
20133 Milano, Italy}
\author[5]{Paola F. Antonietti}
\affil[5]{MOX, Dipartimento di Matematica, Politecnico di Milano, 
Piazza Leonardo da Vinci 32,
20133 Milano, Italy}

\begin{document}
\maketitle

\begin{abstract}
This work aims at presenting a Discontinuous Galerkin (DG) formulation employing a spectral basis for two important models employed in cardiac electrophysiology, namely the monodomain and bidomain models. The use of DG methods is motivated by the characteristic of the mathematical solution of such equations which often corresponds to a highly steep wavefront. Hence, the built-in flexibility of discontinuous methods in developing adaptive approaches, combined with the high-order accuracy, can well represent the underlying physics. The choice of a semi-implicit time integration allows for a fast solution at each time step. The article includes some numerical tests to verify the convergence properties and the physiological behaviour of the numerical solution. Also, a pseudo-realistic simulation turns out to fully reconstruct the propagation of the electric potential, comprising the phases of depolarization and repolarization, by overcoming the typical issues related to the steepness of the wave front.
\end{abstract}

\section{Introduction}

The heart activity is defined through the cardiac cycle which, at a first analysis, is characterized by two alternating phases: the diastole, the period of relaxation, and the systole, the period of contraction. The continuous sequence of contractions, formally called rhythmicity, is caused by the propagation of bio-electrical signals through the cells resulting in an active contraction of the cardiac muscle. Mathematical models are widely employed to understand and predicting the complex processes underlying this phenomenon. The most popular ones are the monodomain and the bidomain models, which mathematically correspond to systems of (non linear) reaction-diffusion equations coupled to a system of ordinary differential equations for the ionic currents \cite{colli_franzone2014,acta}.

The cardiac electrical propagation is represented by a wave characterized by a fast and steep front. Therefore, in order to accurately capture the proper dynamics, Finite Element methods usually require very refined computational grids. This clearly entails high computational costs to solve the electrophysiology problem in real scenarios. On the other hand, high-order spectral methods have become popular for their ability at capturing sharp parts such as shock waves. In this context, Discontinuous Galerkin (DG) methods may provide an effective alternative  for the solution of the electro-physiology problems, guaranteeing high-order accuracy and more flexibility, see \cite{hesthaven2007,saglio2024highorder} and for the problem at hand \cite{hoermann2018adaptive}. 

The objective of this work is to present and test in practice a high-order DG method to discretize the monodomain and bidomain models. In particular,  \Cref{sec:thephysicalproblemanditsmathematicalformulation} is aimed at introducing the electrophysiology problem and the associated mathematical models. The Discontinuous Galerkin formulation is illustrated in \Cref{sec:semi-discretediscontinuousgalerkinapproximtion} while \Cref{sec:anoteonthebasisfunctions} provides the details of the high-order approximation. In order to obtain the final  numerical algorithm, \Cref{sec:fullydiscreteformulation} incorporates the semi-implicit time discretization. Finally, \Cref{sec:numericalresults} presents some numerical tests. Specifically, we perform a preliminary analysis in model problems to highlight the feasibility of this approach showing some convergence tests refining the mesh size and increasing the polynomial order. Then, the results of a realistic (although simplified) simulation of the propagation of the electric potential is shown to highlight the effectiveness of the proposed method. Final remarks are then discussed.

\section{The physical problem and the mathematical formulation}
\label{sec:thephysicalproblemanditsmathematicalformulation}
 In what follows we present a brief introduction to the monodomain and bidomain equations. For a detailed derivation, we address the interested reader to \cite{acta}.

    \subsection{Mathematical models}
    The heart's active mechanical contraction is triggered by the cardiac cell's electrical activation. Cardio-myocytes are activated and deactivated at each heartbeat following a characteristic electrical cycle.
     \noindent The cell is initially at rest ($-90 \, mV$). At the beginning of the activation, its potential increases rapidly ($\approx2\,ms$) and reaches the value of $+20\,mV$. Later, a plateau near $0\,mV$ is observed, followed by a slow repolarization to the initial potential. 
     
    From a microscopical point of view, each single cell is involved in a passage of chemical ions through specific channels, e.g., calcium $Ca^{++}$, sodium $Na^+$ and potassium $K^+$. From a macroscopical point of view, one can describe the dynamics as a continuous electrical diffusion over the entire cardiac tissue driven by the directions of the muscle fibers \cite{piersanti2021}. 
    
    Applying general electromagnetism laws, the bidomain model has been formulated for the macroscopic dynamics (see \cite{colli_franzone2014,acta} for more details). To complete the formulation, a mathematical model for the ionic current is required. In such context, it is noteworthy to mention the original formulation by Hodgkin and Huxley \cite{hodgkinhuxley}. FitzHugh \cite{fitzhugh1961} and Nagumo \cite{nagumo1962} proposed later a simplification of the latter and complex models successively followed such as the ones proposed in \cite{alievpanfilov,ten2004model,ten2006alternans}.
    In this work, the FitzHugh-Nagumo reduced ionic model (FHN) is considered. This simple model represents only one ionic channel o describe the ionic currents. 

    Given an open and bounded domain $\Omega \in \mathbb{R}^d$, $d=2,3$, and a final time $T>0$, the unknowns of the bidomain model are: the trans-membrane potential $V_m = \phi_i - \phi_e$, where $\phi_i$ and $\phi_e$ are the intracellular and extracellular potentials, respectively, and the gating variable $w$ representing the percentage of opening of the ionic channel. The model parameters include: the positive constants $\chi_m,C_m$ , representing the surface area-to-volume ratio and the membrane capacitance, the permeability tensors $\Sigma_i, \Sigma_e$ in the internal and external cellular field, the external applied currents $I_i^{ext}, I_e^{ext}$, and some known constants to tune the ionic model ($\kappa,a,\epsilon,\Gamma$). In particular, $\Sigma_i$ and $\Sigma_e$ account for the anisotropy given by the cardiac fibers. Furthermore, initial and Neumann boundary conditions are imposed through some known functions $\phi_{i,0},\phi_{e,0},w_{0},b_i,b_e$. The former conditions assign the initial state of the system while the latter conditions prescribe the behaviour on the boundary of the domain in terms of inward or outward currents.

     \begin{problem}[Bidomain model coupled with FitzHugh-Nagumo model]\label{prob:bidomain}
      	For each $t\in(0,T]$, find $\phi_i$, $\phi_e$ and $w$ such that:
    	\begin{equation*} 
    	\begin{cases} \displaystyle
    	\chi_m C_m\frac{\partial V_m}{\partial t} - \nabla \cdot (\Sigma_i \nabla \phi_i) + \chi_m I_{ion}(V_m,w) = I_i^{ext}   \qquad & \text{in } \Omega \times (0,T],
    	\\ \displaystyle
    	-\chi_m C_m\frac{\partial V_m}{\partial t} - \nabla \cdot (\Sigma_e \nabla \phi_e) - \chi_m I_{ion}(V_m,w) = -I_e^{ext} \qquad   & \text{in } \Omega \times (0,T],
    	\\ \displaystyle
        V_m = \phi_i - \phi_e & \text{in } \Omega \times [0,T], \\
    	I_{ion}(V_m,w)=\kappa V_m(V_m-a)(V_m-1)+w & \text{in } \Omega \times [0,T],
    	\\\displaystyle
    	\frac{\partial w}{\partial t} = \epsilon(V_m-\Gamma w) \qquad & \text{in } \Omega \times (0,T],
    	\\\displaystyle
     \Sigma_i\nabla \phi_i \cdot \mathbf{n} = b_i \qquad  & \text{on } \partial \Omega \times (0,T],
    	\\\displaystyle
    	\Sigma_e\nabla \phi_e \cdot \mathbf{n} = b_e  \qquad & \text{on } \partial \Omega \times (0,T],
    	\\
    	\phi_i=\phi_{i,0}\qquad & \text{in } \Omega\times\{t=0\}.
     \\
    	\phi_e=\phi_{e,0}\qquad & \text{in } \Omega\times\{t=0\}.
     \\
    	w=w_{0}\qquad & \text{in } \Omega\times\{t=0\}.
    	\end{cases}
    	\end{equation*}
    	\end{problem}
      \noindent Note that Problem \ref{prob:bidomain} needs the following compatibility condition, which is necessary for the existence of the solution:
    \begin{equation} \label{compatibility}
    \int_{\Omega} I_i^{ext} - \int_{\Omega} I_e^{ext} = -\int_{\partial \Omega} b_i - \int_{\partial \Omega} b_e.
    \end{equation}

    A simpler formulation can be derived assuming the intracellular and extracellular permeability tensors to be proportional \cite{colli_franzone2014,acta}. This assumption is valid when the propagation is regular and no chaotic processes such as fibrillation are taking place. In such a way, the unknowns reduce to $V_m$ and $w$. Therefore, the monodomain problem reads as follows:
     \begin{problem}[Monodomain problem coupled with FitzHugh-Nagumo model]\label{prob:monodomain}
      	For each $t\in(0,T]$, find $V_m$ and $w$ such that:
    	\begin{equation*}
    	\begin{cases}\displaystyle
    	\chi_m C_m\frac{\partial V_m}{\partial t} - \nabla \cdot (\Sigma \nabla V_m) + \chi_m I_{ion}(V_m,w) = I^{ext}    & \text{in } \Omega \times (0,T],
    	\\\displaystyle
    	I_{ion}(V_m,w)=\kappa V_m(V_m-a)(V_m-1)+w & \text{in } \Omega \times [0,T],
    	\\\displaystyle
    	\frac{\partial w}{\partial t} = \epsilon(V_m-\Gamma w)  & \text{in } \Omega \times (0,T],
    	\\\displaystyle
     \Sigma\nabla V_m \cdot \mathbf{n} = b \qquad  & \text{on } \partial \Omega \times (0,T],
    	\\
    	V_m=V_{m,0} & \text{in } \Omega\times\{t=0\},
        \\
        w = w_0 & \text{in } \Omega\times\{t=0\},
    	\end{cases}
    	\end{equation*}
    	\end{problem}
\noindent where $\Sigma=\xi/(1+\xi)\Sigma_i, \,
I^{ext}=(\xi I^{ext}_i+I^{ext}_e)(1+\xi)$, $\xi$ being the proportional factor: $\Sigma_e=\xi\Sigma_i$, and for suitable initial and boundary conditions which derive from the bidomain ones. 

\section{Semi-discrete Discontinuous Galerkin approximation}
\label{sec:semi-discretediscontinuousgalerkinapproximtion}
Starting from the strong form given by Problem \ref{prob:bidomain} and Problem \ref{prob:monodomain}, the next step is the pursuit of a Discontinuous Galerkin semi-discrete formulation. Terms and symbols are defined similarly to the usual convention \cite{antonietti2011,arnold2002unified}.

\subsection{Discontinuous Galerkin formulation}
    Let us introduce a shape-regular triangulation $\mathcal{T}_h$ of $\Omega$, where $\mathcal{F} _h=\mathcal{F} _h^I \cup \mathcal{F} _h^B$ is the set of the faces of the partition which includes the internal and boundary faces, respectively. Let the DG space be defined as $\Theta^{h,p} = \{v_h \in L^2(\Omega) : v_h|_\mathcal{K} \in \mathbb{P}^{p}(\mathcal{K})  \; \forall \mathcal{K} \in \mathcal{T}_h \}$, where $\mathbb{P}^p(\mathcal{K})$ is the space of polynomials of total degree less than or equal to $p\ge1$ over $\mathcal{K}\in \mathcal{T}_h$ and $\cdot|_\mathcal{K}$ is the restriction operator to the element $\mathcal{K}$. Moreover, we define $N_h=\dim(\Theta^{h,p})<\infty$ as the dimension of the space $\Theta^{h,p}$ and the following bilinear forms.
    \begin{equation*}\label{forgamma}
    \begin{aligned}\displaystyle
    \langle u_h,v_h\rangle_{\Theta^{h,p}}:=& \sum_{\mathcal{K} \in \mathcal{T}_h} \int_{\mathcal{K}}u_hv_h \; dx,\\
    \newline
    \langle u_h,v_h\rangle_{\mathcal{F} _h^B}:=& \sum_{\mathcal{F} \in \mathcal{F}_h^B} \int_{\mathcal{F}}u_hv_h \; ds,\\
    \newline
    a_{\Sigma_*}(u_h,v_h):=&\sum_{\mathcal{K} \in \mathcal{T}_h} \int_{\mathcal{K}}{\Sigma_*\nabla u_h \cdot \nabla v_h \; dx}-\sum_{F \in \mathcal{F}_h^I} \int_F { \{\!\{\Sigma_* \nabla_h u_h \}\!\} \cdot [\![v_h]\!] \; ds} + \\
    &-\theta \sum_{F \in \mathcal{F}_h^I} \int_F{ \{\!\{\Sigma_* \nabla_h v_h\}\!\} \cdot [\![u_h]\!]\; ds}+\sum_{F \in \mathcal{F}_h^I}\int_F \gamma(\Sigma_*) [\![u_h]\!] \cdot [\![v_h]\!] \; ds.
    \end{aligned}
    \end{equation*}
    Here $\Sigma_*$ is either $\Sigma_i$ or $\Sigma_e$, $\theta \in \{-1,0,1\}$, $\nabla_h$ denotes the element-wise gradient and we set
    
    \begin{equation}
        \label{eq:gamma}\gamma(\Sigma_*) = \overline{\gamma}\,\mathbf{n}^T\Sigma_*\mathbf{n}>0,\qquad \overline{\gamma} := \alpha p^2/h,  
       \end{equation}
    where $\overline{\gamma}$ is a stability parameter defined edge-wise, with $h>0$ the mesh size (supposed to be quasi uniform) and  $\alpha >0$ a fixed parameter,
and $\mathbf{n} \in \mathbb{R}^d$ is the outward normal unitary vector to the corresponding face $F$. Moreover, the jump and average operators $[\![\cdot]\!]$ and $\{\!\{\cdot \}\!\}$ are defined on $F\in \mathcal{F}_h^I$ in the standard way \cite{arnold2002unified}, i.e.,
    \begin{equation*}
        [\![u]\!]:= u_1\mathbf{n}_1 + u_2 \mathbf{n}_2, 
    \end{equation*}
    \begin{equation*}
        \{\!\{\mathbf{v} \}\!\}:=\frac{\mathbf{v}_1+\mathbf{v}_2}{2},
    \end{equation*}
    where the subscript $\{1,2\}$ indicates the evaluation of the variable $u$, the vector variable $\mathbf{v}$ or the normal vector $\mathbf{n}$ with respect to the two adjacent elements $\mathcal{K}_1$ and $\mathcal{K}_2$ such that $\overline{F}=\overline{\mathcal{K}_1}\cap \overline{\mathcal{K}_2}$, $K_1,K_2 \in \mathcal{T}_h$.

    Then, we can write the semi-discretized formulation for the bidomain problem.
    \begin{problem}[Bidomain model - semidiscrete DG formulation]\label{prob:bidomain-dg}
    For any $t\in(0,T]$, find $\phi_i^h(t),\phi_e^h(t),w^h(t) \in \Theta^{h,p}$ such that:
    \begin{equation*}
    \begin{cases}\displaystyle
    \left\langle\chi_m C_m \frac{\partial V^h_m}{\partial t}, v_h \right\rangle_{\Theta^{h,p}} +a_{\Sigma_i}\left(\phi_i^h,v_h\right)+
    \left\langle I_{ion}\left(V_m^h,w^h\right),v_h\right\rangle_{\Theta^{h,p}}=\left\langle I_i^{ext}, v_h \right\rangle_{\Theta^{h,p}}+ \left\langle b_i, v_h\right\rangle_{\mathcal{F} _h^B}\\\displaystyle
    \left\langle\chi_m C_m \frac{\partial V^h_m}{\partial t}, v_h \right\rangle_{\Theta^{h,p}} -a_{\Sigma_e}\left(\phi_e^h,v_h\right)+
    \left\langle I_{ion}\left(V_m^h,w^h\right),v_h\right\rangle_{\Theta^{h,p}}=\left\langle I_e^{ext}, v_h \right\rangle_{\Theta^{h,p}}- \langle b_e, v_h\rangle_{\mathcal{F} _h^B}\\\displaystyle
    \displaystyle \left\langle \frac{\partial w^h}{\partial t}, v_h \right\rangle_{\Theta^{h,p}}= \left\langle \epsilon \left(V^h_m - \Gamma w^h\right), v_h \right\rangle_{\Theta^{h,p}}, \\\displaystyle
    \phi_i^h(0) = \phi_{i,0}^h, \\\displaystyle
    \phi_e^h(0) = \phi_{e,0}^h, \\ \displaystyle
    w^h(0) = w^h_0, \\
    \end{cases}
    \end{equation*}
    for each $v_h\in \Theta^{h,p}$,  where $V^h_m=\phi_i^h-\phi_e^h$, $I_{ion}(\cdot,\cdot)$ is defined in Problem \ref{prob:bidomain}, $\phi_{i,0}^h,\phi_{e,0}^h,w^h_0 \in \Theta^{h,p}$ are suitable projections onto $\Theta^{h,p}$ of $\phi_{i,0},\phi_{e,0},w_0$.
    \end{problem}

     \noindent Following the same notation, we also present the weak formulation for the monodomain problem.
    \begin{problem}[Monodomain model - DG semidiscrete formulation]\label{prob:monodomain-dg}
    For any $t\in(0,T]$, find $V_m^h(t) \in \Theta^{h,p}$  and  $w^h(t) \in \Theta^{h,p}$ such that:
    \begin{equation*}
\begin{cases}\displaystyle
\left\langle\chi_m C_m \frac{\partial V^h_m}{\partial t}, v_h \right\rangle_{\Theta^{h,p}} +  a_\Sigma\left(V_m^h,v_h\right) + \left\langle I_{ion}\left(V_m^h,w^h\right),v_h\right\rangle_{\Theta^{h,p}} = \left\langle I^{ext}, v_h \right\rangle_{\Theta^{h,p}}+ \langle b, v_h\rangle_{\mathcal{F} _h^B}, \\\displaystyle
\displaystyle \left\langle \frac{\partial w^h}{\partial t}, v_h \right\rangle_{\Theta^{h,p}}= \left\langle \epsilon \left(V^h_m - \Gamma w^h\right), v_h \right\rangle_{\Theta^{h,p}}, \\\displaystyle
V_m^h(0) = V_{m,0}^h, \\\displaystyle
w^h(0) = w^h_0, 
\end{cases}
\end{equation*}
for every $v_h \in \Theta^{h,p}$, where $I_{ion}$ is defined in Problem \ref{prob:monodomain} and $V_{m,0}^h, w^h_0 \in \Theta^{h,p}$ are suitable projections onto $\Theta^{h,p}$ of $V_{m,0},w_0$.
    \end{problem}
    \noindent In Problem \ref{prob:bidomain-dg} and Problem \ref{prob:monodomain-dg}, three different methods are employed according to the choice of the coefficient $\theta$:
    \begin{itemize}
    \item $\theta=1$: Symmetric Interior Penalty method (SIP) \cite{douglas2008interior};
    \item $\theta=0$: Incomplete Interior Penalty method (IIP) \cite{sun2005symmetric};
    \item $\theta=-1$: Non Symmetric Interior Penalty method (NIP) \cite{riviere2001priori}.
    \end{itemize}
    \noindent Notice that $\alpha$ in the definition \eqref{eq:gamma} needs to be large enough to guarantee coercivity of the SIP and IIP formulations.

 \subsection{Algebraic formulation}
 Let $\{\varphi_j\}_{j=1}^{N_h}$ be a basis of $\Theta^{h,p}$ and let $\boldsymbol{V}^h_m(t), \boldsymbol{\phi}^h_i(t), \boldsymbol{\phi}^h_e(t), \boldsymbol{w}^h(t) \in \mathbb{R}^{N_h}$ be the vectors containing the expansion coefficients of $V_m^h(t), \phi_i^h(t), \phi_e^h(t), w^h(t)$ with respect to such a basis for each instant of time $t\in (0,T]$. 

For $i,j=1,\cdots,N_h$ and $z=\{i,e,\emptyset\}$, we define the following matrices and right hand side vector:
 {\allowdisplaybreaks
 \begin{align*}
\relax[K_z]_{jk} &= \sum_{\mathcal{K} \in \mathcal{T}_h}\int_{\mathcal{K}}\nabla_h\varphi_k \cdot \Sigma_z \nabla_h \varphi_j,\\ 
[W_z]_{jk} &= \sum_{F \in \mathcal{F}_h^I} \int_{F} [\![\varphi_k]\!] \cdot \{\!\{\Sigma_k \nabla_h \varphi_j\}\!\},\\
[S_z]_{jk} &= \sum_{F \in \mathcal{F}_h^I} \int_{F}\gamma(\Sigma_z)[\![\varphi_k]\!] \cdot [\![\varphi_j]\!],\\
[A_z]_{jk} &= [K_z]_{jk}-[W_z]_{kj}-\theta [W_z]_{jk} +[S_z]_{jk} \hspace{7.5mm}\text{stiffness matrix},\\
[M]_{jk} &= \sum_{\mathcal{K} \in \mathcal{T}_h}\int_{\mathcal{K}}\varphi_k\varphi_j \hspace{41mm}\text{mass matrix},\\
[C(V^h_m)]_{jk} &=  \sum_{\mathcal{K} \in \mathcal{T}_h} \int_{\mathcal{K}} \chi_m \kappa(V^h_m-1)(V^h_m-a)\varphi_k\varphi_j \qquad{\text{non-linear reaction matrix}},\\
[\mathbf{R_z}]_j &= \sum_{\mathcal{K} \in \mathcal{T}_h}\int_{\mathcal{K}} I_z^{ext}\varphi_j + \sum_{F \in \mathcal{F}_h^B} \int_F b_z\varphi_j  \hspace{15mm}\text{forcing term}.
\end{align*}}
\vspace{3mm} \\
The semi-discrete DG algebraic formulation of the bidomain problem leads to the following ODE system:
 \begin{problem}[Algebraic formulation of the bidomain model] \label{block_matrix}
 Find $\boldsymbol{V}^h_m(t)= \boldsymbol{\phi}^h_i(t)- \boldsymbol{\phi}^h_e(t), \boldsymbol{w}^h(t) \in \mathbb{R}^{N_h}$ such that for any $t \in (0,T]$:
 \begin{equation*}
 \begin{cases}

 \chi_mC_m \begin{bmatrix}M &-M \\ -M & M \end{bmatrix}
	\begin{bmatrix}\dot{\boldsymbol{\phi}_i^h}(t) \\ \dot{\boldsymbol{\phi}_e^h}(t) \end{bmatrix} 
	 +
	   \begin{bmatrix}A_i + C(V_m^h(t)) & -C(V_m^h(t)) \\ -C(V_m^h(t)) & A_e+ C(V_m^h(t)) \end{bmatrix} 
	   \begin{bmatrix} \boldsymbol{\phi}_i^h(t) \\ \boldsymbol{\phi}_e^h(t)  \end{bmatrix} 
	   +\chi_m \begin{bmatrix}M & 0 \\ 0 & -M \end{bmatrix} 
	   	\begin{bmatrix}\boldsymbol{w}^h(t) \\ \boldsymbol{w}^h(t) \end{bmatrix} = 
	   	\begin{bmatrix} \mathbf{R_i}(t) \\ \mathbf{R_e}(t)\end{bmatrix},
	   	\\
	  M \dot{\boldsymbol{w}}^h(t)= \epsilon M (\boldsymbol{V}_m^h(t)-\Gamma \boldsymbol{w}^h(t)),
\end{cases}
\end{equation*}
together with initial conditions. 
\end{problem}
\noindent The monodomain DG ODE system can be written similarly.
\begin{problem}[Algebraic formulation of the monodomain model]\label{algebraic}
Find $\boldsymbol{V}_m^h(t),\boldsymbol{w}^h(t) \in \mathbb{R}^{N_h}$ such that for any $t \in (0,T]$:
\begin{equation*}
\begin{cases}
\chi_m C_m M \dot{\boldsymbol{V}_m^h}(t) + A \boldsymbol{V}_m^h(t) + C(V_m^h(t)) \boldsymbol{V}_m^h(t) + \chi_m M \boldsymbol{w}^h(t) = \mathbf{R}(t),\\
M \dot{\boldsymbol{w}^h}(t) = \epsilon M (\boldsymbol{V}_m^h(t) - \Gamma \boldsymbol{w}^h(t)), \\
\end{cases}
\end{equation*}
together with initial conditions.
\end{problem}

\section{A note on the basis functions}
\label{sec:anoteonthebasisfunctions}

   The Dubiner spectral functions \cite{dubiner91} are chosen as basis $\{\varphi_j\}_{j=1}^{N_h}$ for the space $\Theta^{h,p}$. The advantage of using an $L^2$ orthogonal basis on simplices rather than a nodal basis relies on the feasibility of applying high-order discretizations without incurring in high costs constructing the discretized space and/or ill-conditioned matrices. Moreover, the mass matrix is in this case diagonal. 

   \begin{figure}[h!]
\centering
\begin{tikzpicture}
\node[inner sep=0pt] at (0,0) {\includegraphics[width=.25\textwidth]{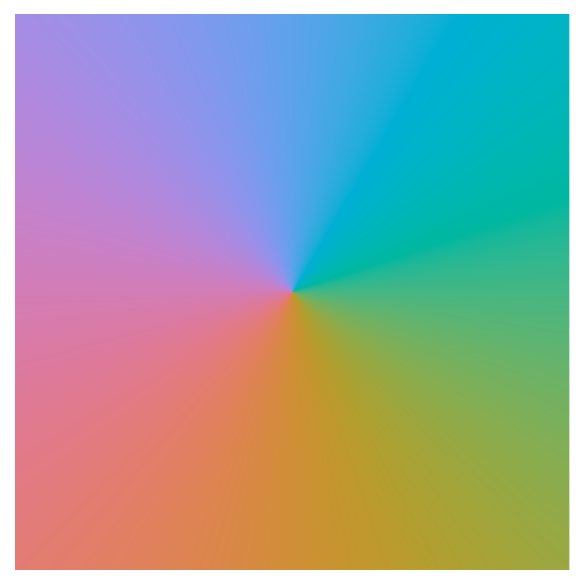}};
\node[inner sep=0pt] at (0.55\textwidth,0) {\includegraphics[width=.25\textwidth]{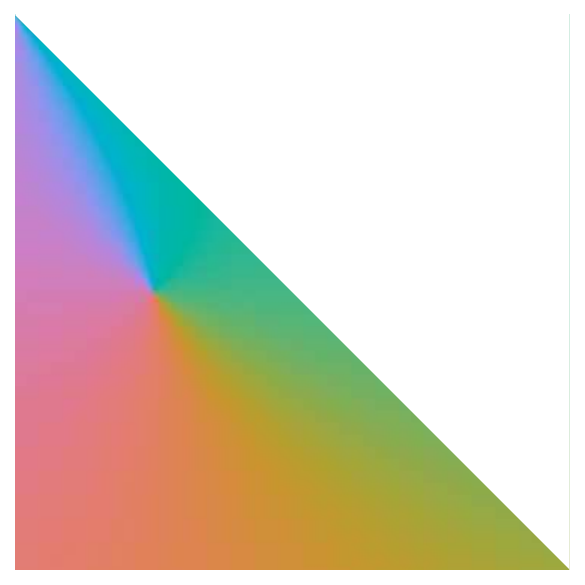}};

\draw[->] (0.15\textwidth,0) -- (0.37\textwidth,0) node[above,midway] {$\hat{Q}\rightarrow\hat{K}$};
\draw[->] (0.15\textwidth,0) -- (0.37\textwidth,0) node[below,midway] {$(a,b)\rightarrow(\xi,\eta)$};

\node[] at (-0.12\textwidth,-0.13\textwidth) {$-1$};
\node[] at (+0.125\textwidth,-0.13\textwidth) {$1$};
\node[] at (-0.135\textwidth,-0.115\textwidth) {$-1$};
\node[] at (-0.13\textwidth,+0.125\textwidth) {$1$};
\node[] at (0,-0.13\textwidth) {$a$};
\node[] at (-0.13\textwidth,0) {$b$};
\node[] at (0.42\textwidth,+0.13\textwidth) {$1$};
\node[] at (0.42\textwidth,-0.11\textwidth) {$0$};
\node[] at (0.44\textwidth,-0.13\textwidth) {$0$};
\node[] at (0.675\textwidth,-0.13\textwidth) {$1$};
\node[] at (0.55\textwidth,-0.13\textwidth) {$\xi$};
\node[] at (0.42\textwidth,0) {$\eta$};

\end{tikzpicture}
\caption{Color plot of the transformation (\ref{transformation_formula}) underlying the construction of the 2D Dubiner spectral basis on simplices. The collapse of the square can be seen on the upper edge of the triangle.}
\end{figure}
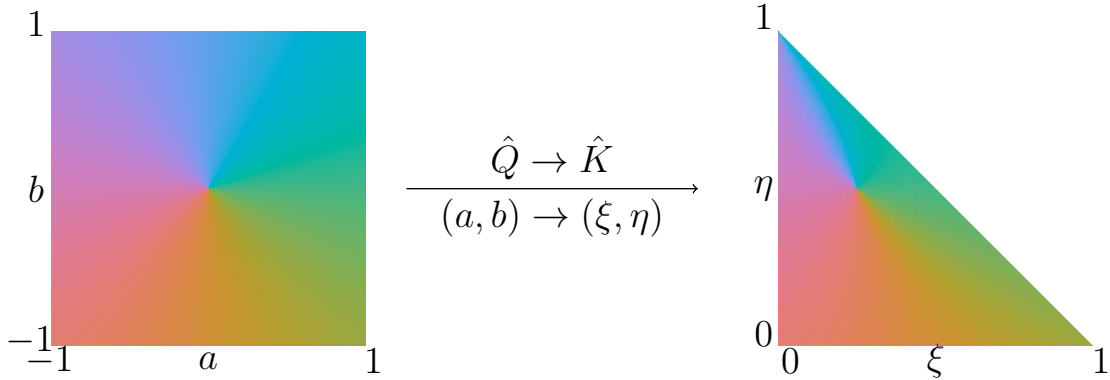

On the other hand, the projection of known functions on the spectral space is certainly more involved since it requires $L^2$ scalar products. This issue can be overcome by implementing fast evaluations of Jacobi polynomials and efficient quadrature formulae to be also used in the weak formulations of Problems \ref{prob:bidomain-dg} and \ref{prob:monodomain-dg}.  

    In $d=2$, the Dubiner basis is defined on the reference simplex as follows \cite{antonietti2011}. This is possible by collapsing the two-dimensional reference square $\hat{Q}=\{ (a, b) : \, -1 \le a \le 1, \, -1 \le b \le 1 \}$ into the reference triangle $\hat{K}=\{ (\xi, \eta) : \, \xi, \eta \ge 0, \,	\xi+\eta \le 1 \} $ applying the following transformation:
    \begin{equation}\label{transformation_formula}
    \xi=\frac{(1+a)(1-b)}{4},  \quad \eta=\frac{(1+b)}{2}.
    \end{equation}
    
    Hence, the Dubiner basis functions are constructed as the transformations of suitable basis functions initially defined on $\hat{Q}$. More precisely, on the reference square $\hat{Q}$ the basis functions are defined as modified tensor products of Jacobi polynomials. 
    
    \begin{definition}[Jacobi polynomials \cite{jacobi}]\label{jacobi}
    The $n$-th Jacobi polynomial of indices $\alpha,\beta \ge -1$ is defined as:
    \begin{equation*}
        P_n^{\alpha,\beta}(x)= \frac{(-1)^n}{2^n n!} (1-x)^{-\alpha}(1+x)^{-\beta} \frac{d^n}{dx^n} \left[(1-x)^\alpha(1+z)^\beta(1-z^2)^n\right], n\ge0.
    \end{equation*}
    \end{definition}
    \begin{definition}[2D Dubiner basis functions \cite{antonietti2011}]
    If $d=2$, the Dubiner basis function indexed by $(i,j) \in \mathbb{N}^2, i+j\le p$ is defined as:
    \begin{gather*}
        \varphi_{i,j}:\hat{K}\rightarrow\mathbb{R}\\
        \varphi_{i,j}(\xi,\eta)=c_{i,j} (1-\eta)^i P_i^{0,0}\left(\frac{2\xi}{1-\eta} -1\right) P_j^{2i+1,0}(2\eta-1),
    \end{gather*}
    where $c_{i,j}=\sqrt{2(2i+1)(i+j+1)}$.
    \end{definition}

\begin{remark} 
Despite the great potentiality of an orthogonal spectral basis, it is often necessary to pass to nodal representations. Let $\{\varphi_j\}_{j=1}^{N_h}$ be the set of Dubiner basis functions and let $v_h=\sum_{j=1}^{N_h}v_j\varphi_j \in \Theta^{h,p}$. A nodal evaluation of such solution can be performed through the linear combination
	\begin{equation} \label{ref3}
	v_h(x) = \sum_{j=1}^{N_h} v_j\varphi_j(x), \: \forall x \in \Omega
	\end{equation}
	while the converse operation can be performed through a $L^2$ scalar product thanks to the orthonormality property:
 \begin{equation}\label{ref4}
	\begin{gathered}
	v_j = \int_\Omega v_h(x) \varphi_j(x) \,dx \hspace{4mm} \forall j=1,\dots,N_h.
	\end{gathered}
	\end{equation}
 \end{remark}

 \section{Fully discrete formulation}
 \label{sec:fullydiscreteformulation}
In this section a first order semi-implicit time discretization is presented in order to obtain a fully discretized system of equations from Problem \ref{prob:bidomain-dg} and \ref{prob:monodomain-dg}. Thus, we split the interval $(0,T]$ into $N$ uniform sub-intervals $(t^n,t^{n+1}]$ of length $\Delta t$ such that $t^n=n \Delta t \quad \forall n=0,\cdots,N-1$. Then, we define the fully discretized solutions $V_m^{h,n},w^{h,n},\phi_i^{h,n}, \phi_e^{h,n} \in \Theta^{h,p}$ that are approximations of $V_m^h(t^n)$, $w^h(t^n)$, $\phi_i^h(t^n)$, $\phi_e^h(t^n)$, respectively, whose expansion coefficients vectors are denoted with $\boldsymbol{V}_m^{h,n}$, $\boldsymbol w^{h,n}$, $\boldsymbol{\phi}_i^{h,n}$,$\boldsymbol{\phi}_e^{h,n} \in \mathbb{R}^{N_h}$, respectively. 

In the semi-implicit scheme, only the non-linear term is treated in explicit in order to take advantage of an implicit discretization while still having linearity in the time-step advancing scheme. 
 \begin{problem}[Fully discretized formulation of the Bidomain model] \label{prob:bidomain-full} Given  $\phi_i^{h,0},\phi_e^{h,0},w^{h,0} \in \Theta^{h,p}$, find $\boldsymbol{V}_m^{h,n+1} = \boldsymbol{\phi}_i^{h,n+1} -\boldsymbol{\phi}_e^{h,n+1},\boldsymbol{w}^{h,n+1} \in \mathbb{R}^{N_h}$ for $n=0, \ldots,N-1$ such that:
 \begin{equation*}
 \begin{cases}
 \begin{aligned}
 \displaystyle\frac{\chi_mC_m}{\Delta t} \begin{bmatrix}M &-M \\ -M & M \end{bmatrix}
	\begin{bmatrix}\boldsymbol{\phi}_i^{h,n+1}-\boldsymbol{\phi}_i^{h,n} \\ {\boldsymbol{\phi}_e^{h,n+1}-\boldsymbol{\phi}_e^{h,n}} \end{bmatrix}
	 + \begin{bmatrix}A_i + C(V_m^{h,n}) & -C(V_m^{h,n}) \\ -C(V_m^{h,n}) & A_e + C(V_m^{h,n}) \end{bmatrix} 
	 \begin{bmatrix}\boldsymbol{\phi}_i^{h,n+1} \\ \boldsymbol{\phi}_e^{h,n+1} \end{bmatrix} & \\
	   +\chi_m \begin{bmatrix}M & 0 \\ 0 & -M \end{bmatrix} 
	   	\begin{bmatrix}\boldsymbol{w}^{h,n+1} \\ \boldsymbol{w}^{h,n+1} \end{bmatrix} = 
	   	\begin{bmatrix} \mathbf{R_i}(t^{n+1}) \\ \mathbf{R_e}(t^{n+1})\end{bmatrix}&,
	   	\end{aligned}\\
	   	\vspace{3mm} 
	  M \mathlarger{\frac{\boldsymbol{w}^{h,n+1}-\boldsymbol{w}^{h,n}}{\Delta t}}= \epsilon M (\boldsymbol{V}_m^{h,n}-\Gamma \boldsymbol{w}^{h,n+1}).
\end{cases}
\end{equation*}
\end{problem}
\begin{problem}[Fully discretized formulation of the Monodomain model]
Given $V_m^{h,0},w^{h,0} \in \Theta^{h,p}$, find $\boldsymbol{V}_m^{h,n+1},\boldsymbol{w}^{h,n+1} \in \mathbb{R}^{N_h}$ for $n=0, \ldots, N-1$ such that:
\begin{equation*}
\begin{cases}\displaystyle
\chi_m C_m M \frac{\boldsymbol{V}_m^{h,n+1}-\boldsymbol{V}_m^{h,n}}{\Delta t} + A \boldsymbol{V}_m^{h,n+1} + C(V_m^{h,n}) \boldsymbol{V}_m^{h,n+1} + \chi_m M \boldsymbol{w}^{h,n+1} = \mathbf{R}(t^{n+1}),\\\displaystyle
M \frac{\boldsymbol{w}^{h,n+1}-\boldsymbol{w}^{h,n}}{\Delta t} = \epsilon M (\boldsymbol{V}_m^{h,n} - \Gamma \boldsymbol{w}^{h,n+1}).
\end{cases}
\end{equation*}
\end{problem}

\begin{remark}
Problem \ref{prob:bidomain} and \ref{prob:monodomain} are known to be well-posed in the sense of existence, uniqueness and regularity of the solutions under suitable assumptions \cite{colli_franzone2014,bourgault2009}. More precisely, $V_m(t)$ is unique in $H^1(\Omega)$ $\forall t \in (0,T]$ while $\phi_i(t)$ and $\phi_e(t)$ are unique in $H^1(\Omega)/\mathbb{R}$, i.e. up to an additive constant. Since $V_m$ is defined as the difference between the two potentials, this constant is necessarily the same for both $\phi_i$ and $\phi_e$. Therefore, the numerical solution of Problem \ref{prob:bidomain-full} requires to fix the constant, e.g., by imposing the value of $\phi_i$ or $\phi_e$ in a particular point of the domain or by imposing its mean value. In the former case, the condition
\begin{equation*}
    \phi_i(\mathbf{x},t^n) = c \hspace{3mm} \text{or} \hspace{3mm} \phi_e(\mathbf{x},t^n) = c
\end{equation*}
is applied for arbitrary $\mathbf{x}\in\Omega$, $c\in\mathbb{R}$, $\forall n=1,\dots,N$. We opt instead for the latter strategy which is more consistent with the variational formulation of the problem in $H^1(\Omega)$ and is represented by the condition 
\begin{equation*}
    \sum_{j=1}^{N_h} \boldsymbol{\phi}_{i,j}^{h,n} d_j = c \hspace{3mm} \text{or} \hspace{3mm} \sum_{j=1}^{N_h} \boldsymbol{\phi}_{e,j}^{h,n} d_j = c,
\end{equation*}
where $d_j:=\int_\Omega \varphi_j$, for an arbitrary $c\in\mathbb{R}$, $\forall n=1,\dots,N$.
\end{remark}

 \section{Numerical results}
  \label{sec:numericalresults}
This section is devoted to some numerical tests aimed at assessing the quality of the proposed numerical strategies. Some convergence tests are discussed in \Cref{sec:verificationtestcases} while a more realistic simulation is shown in \Cref{sec:towardsarealisticsimulation}. The code has been implemented in a small C++ library, called \texttt{DUBeat}\footnote{\url{https://matteocalafa.com/DUBeat/}}, which is largely based on \texttt{lifex} \cite{africa2022}.

 \subsection{Verification test cases}

 \begin{figure}[h!]
    \centering
    \includegraphics[width=0.5\textwidth]{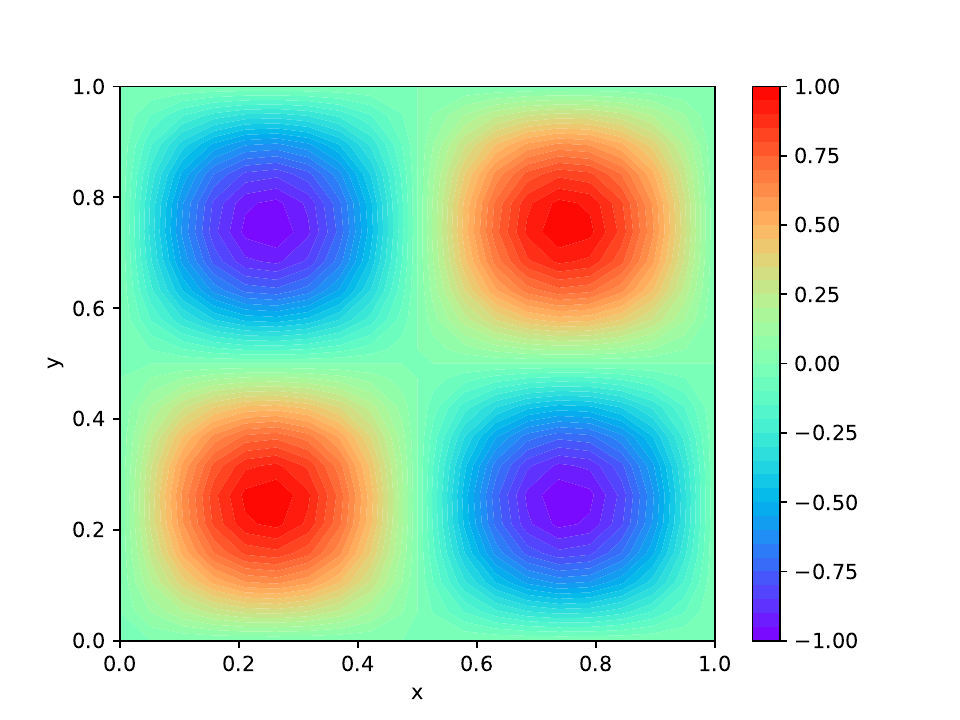}
    \caption{Contour plot of the exact solution for $\phi_e$ in \Cref{sec:verificationtestcases}.}
\end{figure}

\label{sec:verificationtestcases}
We consider a simple square domain $\Omega=(0,1)^2 \subset \mathbb{R}^2$. We opt for the symmetric interior penalty method (SIP, $\theta=1$) and the following parameters choices: $T=3 \cdot 10^{-3}$, $\Delta T = 10^{-4}$, $\chi_m = 10^5$, $C_m=1$, $\kappa = 19.5$, $\epsilon=1.2$, $\Gamma=0.1$, $a=1.3\cdot 10^{-2}$, $\alpha=10$. Moreover, $\Sigma = \Sigma_i = \Sigma_e = 0.12 \mathcal{I}_2$, where $\mathcal{I}_2 \in M_{2,2}(\mathbb{R})$ is the $2$-dimensional identity matrix, i.e. we assume that the coordinate system is already aligned with the principal fibers directions and there is no anisotropy. 

Boundary conditions and applied currents are assigned assuming that the exact solutions are
\begin{equation*}
    \begin{aligned}
        \phi_i (x,y,t) = 2\sin(2\pi x)\sin(2\pi y)e^{-5t}, \hspace{3mm} \phi_e (x,y,t) = \sin(2\pi x)\sin(2\pi y)e^{-5t}.
    \end{aligned}
\end{equation*}

In the following tests, the errors are computed with respect to the standard $L^2(\Omega)$ norms, the $H^1(\Omega)$ norm $\| v\|^2_{H^1(\Omega)}:= \| v\|^2_{L^2(\Omega)} + \| |\nabla v| \|^2_{L^2(\Omega)}$ and the DG norm
\begin{equation*}
    \|v\|^2_{DG}:= \|\nabla_h v \|_{L^2(\Omega)} + \| \overline{\gamma}^{1/2} [\![v]\!]\|_{L^2(\mathcal{F}_h)}^2.
\end{equation*}
 
\begin{figure}[ht]
    \centering
    
\begin{tikzpicture}[scale=2.7]
\draw (0,0) -- (0,1);
\draw (0,1) --  (1,1);
\draw (1,1) -- (1,0);
\draw (1,0) -- (0,0);
\draw (0,0) -- (1,1);
\draw [->] (1.2,0.5) -- (1.5,0.5);

\draw (1.8,0) -- (1.8,1);
\draw (1.8,1) --  (2.8,1);
\draw (2.8,1) -- (2.8,0);
\draw (2.8,0) -- (1.8,0);
\draw (1.8,0) -- (2.8,1);
\draw (2.3,0) -- (2.3,0.5);
\draw (2.3,0.5) -- (2.8,0.5);
\draw (2.3,0) -- (2.8,0.5);
\draw (1.8,0.5) -- (2.3,0.5);
\draw (2.3,0.5) -- (2.3,1);
\draw (1.8,0.5) -- (2.3,1);
\draw [->] (3,0.5) -- (3.3,0.5);

\draw (3.5,0) -- (3.5,1);
\draw (3.5,1) --  (4.5,1);
\draw (4.5,1) -- (4.5,0);
\draw (4.5,0) -- (3.5,0);
\draw (3.5,0) -- (4.5,1);
\draw (3.5,0.75) -- (3.75,1);
\draw (3.5,0.5) -- (4,1);
\draw (3.5,0.25) -- (4.25,1);
\draw (3.75,0) -- (4.5,0.75);
\draw (4,0) -- (4.5, 0.5);
\draw (4.25,0) -- (4.5, 0.25);
\draw (3.75,0) -- (3.75,1);
\draw (4,0) -- (4,1);
\draw (4.25,0) -- (4.25,1);
\draw (3.5,0.25) -- (4.5,0.25);
\draw (3.5,0.5) -- (4.5,0.5);
\draw (3.5,0.75) -- (4.5,0.75);

\fill[gray, opacity=0.2] (0,0) -- (0,1) -- (1,1) -- (1,0) -- (0,0);
\fill[gray, opacity=0.2] (1.8,0) -- (1.8,1) -- (2.8,1) -- (2.8,0) -- (1.8,0);
\fill[gray, opacity=0.2] (3.5,0) -- (3.5,1) -- (4.5,1) -- (4.5,0) -- (3.5,0);
\end{tikzpicture}
\caption{Sample of uniformly refined grids with granularity $h\approx 2^{-\sigma}$, $\sigma=0,1,2$.}

\label{fig:refinement}
\end{figure}
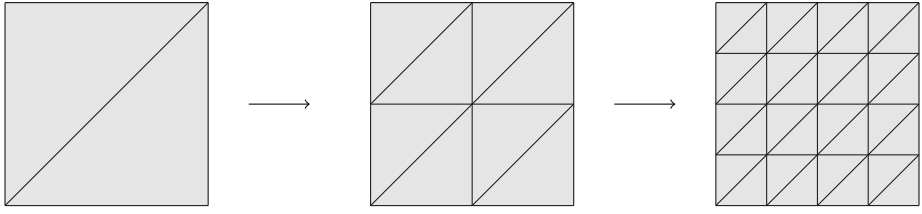

 We tested the proposed formulations in two dimensions on a sequence of uniformly refined grids of granularity $h\approx2^{-\sigma}$, $\sigma=0,1,\dots,5$, cf. \Cref{fig:refinement}.

\begin{figure}
\centering 
\begin{tikzpicture}[scale=0.80]
\begin{semilogyaxis}[
    xlabel={Polynomial degree $p$},
    ylabel={Error},
    xmin=1, xmax=5,
    ymin=0.0001, ymax=6,
    xtick={1,2,3,4,5},
    legend pos=south west,
    ymajorgrids=true,
    grid style=dashed,
]

\addplot[
    color=blue,
    mark=*,
    ]
    coordinates {
    (1,0.201083) (2,0.261212) (3,0.034069) (4,0.007991) (5,0.001307)
    };

\addplot[
    color=red,
    mark=x,
    ]
    coordinates {
    (1,0.089865) (2,0.037163) (3,0.003939) (4,0.000888) (5,0.000252)
    };

\addplot[
    color=green,
    mark=square*,
    ]
    coordinates {
    (1,2.349479) (2,1.273132) (3,0.233224) (4,0.061651) (5,0.009736)
    };

\addplot[
    color=cyan,
    mark=triangle*,
    ]
    coordinates {
    (1,4.201246) (2,3.155599) (3,0.581762) (4,0.160401) (5,0.025470)
    };
    \legend{$L^\infty$, $L^2$, $H^1$, $DG$}
    
\end{semilogyaxis}
\end{tikzpicture}
    \caption{Computed errors vs polynomial degree $p$ (semilogy scale) for monodomain model, $h\approx 2^{-3}$.}
    \label{fig:conv_p}
\end{figure}
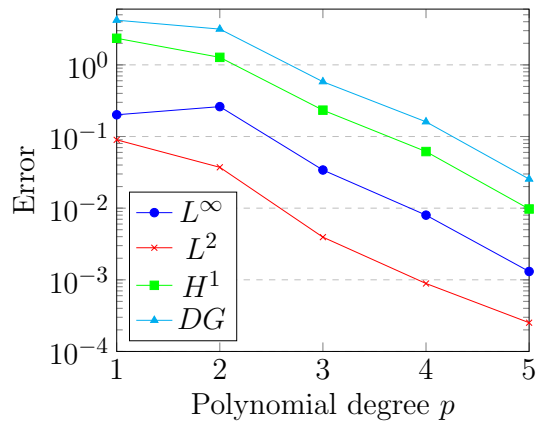

\Cref{fig:mono21,fig:mono22,fig:bido21,fig:bido22} show the computed errors obtained with the monodomain and bidomain problems in two dimensions. Except for the pre-asymptotic regime when $h\approx 1$, the errors for $d=2$ shows the expected theoretical convergence orders that are shown in the bottom triangles. If $p\in \mathbb{N}$ is the polynomial degree, the order is indeed $O(h^p)$ for the $H^1$ and $DG$ errors and $O(h^{p+1})$ for the $L^\infty$ and $L^2$ errors. We also observe that the curves are more flat for $h\ll1,p=2$ because the error due to the grid refinement is very small and therefore comparable to the error caused by the time advancing scheme. \\
Correct convergence orders are also verified for $p$-refinements in \cref{fig:conv_p} where the grid size $h=2^{-3}$ is fixed and the polynomial degree varies from 1 to 5.

\begin{figure}

 \begin{subfigure}{0.5\textwidth}
 \centering
 \begin{tikzpicture}[scale=0.80]
\begin{loglogaxis}[
    xlabel={Grid length ($h$)},
    ylabel={Error},
    xmin=0.03125, xmax=1,
    ymin=0.001, ymax=6,
    xtick={0.03125,0.0625,0.125,0.25,0.5,1},
    legend pos=south east,
    ymajorgrids=true,
    grid style=dashed,
]

\addplot[
    color=blue,
    mark=*,
    ]
    coordinates {
    (1,0.435193) (0.5,0.478846) (0.25,0.201083) (0.125,0.057323) (0.0625,0.012172) (0.03125,0.001973)
    };

\addplot[
    color=red,
    mark=x,
    ]
    coordinates {
    (1,0.248508) (0.5,0.329477) (0.25,0.089865) (0.125,0.02206) (0.0625,0.005203) (0.03125,0.001269)
    };

\addplot[
    color=green,
    mark=square*,
    ]
    coordinates {
    (1,4.347388) (0.5,3.453478) (0.25,2.349479) (0.125,1.277812) (0.0625,0.632654) (0.03125,0.319929)
    };

\addplot[
    color=cyan,
    mark=triangle*,
    ]
    coordinates {
    (1,4.881491) (0.5,4.38083) (0.25,4.201246) (0.125,1.61246) (0.0625,0.745052) (0.03125,0.337791)
    };
    \legend{$L^\infty$, $L^2$, $H^1$, $DG$}
    
 \logLogSlopeTriangle{0.6}{0.1}{0.25}{1}{black};
 \logLogSlopeTriangle{0.6}{0.1}{0.1}{2}{black};
    
\end{loglogaxis}
\end{tikzpicture}
\caption{Monodomain model, $p=1$}
\label{fig:mono21}
\end{subfigure}
\begin{subfigure}{0.5\textwidth}
\centering 
\begin{tikzpicture}[scale=0.80]
\begin{loglogaxis}[
    xlabel={Grid length ($h$)},
    ylabel={Error},
    xmin=0.03125, xmax=1,
    ymin=0.0001, ymax=15,
    xtick={0.03125,0.0625,0.125,0.25,0.5,1},
    legend pos=south east,
    ymajorgrids=true,
    grid style=dashed,
]

\addplot[
    color=blue,
    mark=*,
    ]
    coordinates {
    (1,1.531658) (0.5,0.638646) (0.25,0.261212) (0.125,0.017634) (0.0625,0.002271) (0.03125,0.000594)
    };

\addplot[
    color=red,
    mark=x,
    ]
    coordinates {
    (1,0.517353) (0.5,0.147708) (0.25,0.037163) (0.125,0.002959) (0.0625,0.000401) (0.03125,0.000248)
    };

\addplot[
    color=green,
    mark=square*,
    ]
    coordinates {
    (1,4.820600) (0.5,2.515431) (0.25,1.273132) (0.125,0.249169) (0.0625,0.057774) (0.03125,0.014409)
    };

\addplot[
    color=cyan,
    mark=triangle*,
    ]
    coordinates {
    (1,11.359070) (0.5,7.8747) (0.25,3.155599) (0.125,0.507768) (0.0625,0.089558) (0.03125,0.024607)
    };
    \legend{$L^\infty$, $L^2$, $H^1$, $DG$}
    
 \logLogSlopeTriangle{0.6}{0.1}{0.25}{2}{black};
 \logLogSlopeTriangle{0.6}{0.1}{0.1}{3}{black};
    
\end{loglogaxis}
\end{tikzpicture}
    \caption{Monodomain model, $p=2$}
    \label{fig:mono22}
\end{subfigure}
\\

 \begin{subfigure}{0.5\textwidth}
 \centering
 \begin{tikzpicture}[scale=0.80]
\begin{loglogaxis}[
    xlabel={Grid length ($h$)},
    ylabel={Error},
    xmin=0.03125, xmax=1,
    ymin=0.001, ymax=6,
    xtick={0.03125,0.0625,0.125,0.25,0.5,1},
    legend pos=south east,
    ymajorgrids=true,
    grid style=dashed,
]

\addplot[
    color=blue,
    mark=*,
    ]
    coordinates {
    (1,0.7314) (0.5,0.4564) (0.25,0.4682) (0.125,0.1467) (0.0625,0.0385) (0.03125,0.0094)
    };

\addplot[
    color=red,
    mark=x,
    ]
    coordinates {
    (1,0.5274) (0.5,0.2664) (0.25,0.1166) (0.125,0.0334) (0.0625,0.0084) (0.03125,0.0021)
    };

\addplot[
    color=green,
    mark=square*,
    ]
    coordinates {
    (1,4.2411) (0.5,0.37775) (0.25,2.6411) (0.125,1.3911) (0.0625,0.6754) (0.03125,0.3213)
    };

\addplot[
    color=cyan,
    mark=triangle*,
    ]
    coordinates {
    (1,4.6149) (0.5,4.2403) (0.25,4.4222) (0.125,2.1106) (0.0625,1.0303) (0.03125,0.4838)
    };
    \legend{$L^\infty$, $L^2$, $H^1$, $DG$}
    
 \logLogSlopeTriangle{0.6}{0.1}{0.25}{1}{black};
 \logLogSlopeTriangle{0.6}{0.1}{0.1}{2}{black};
    
\end{loglogaxis}
\end{tikzpicture}
\caption{Bidomain model, $p=1$}
\label{fig:bido21}
\end{subfigure}
\begin{subfigure}{0.5\textwidth}
\centering 
\begin{tikzpicture}[scale=0.80]
\begin{loglogaxis}[
    xlabel={Grid length ($h$)},
    ylabel={Error},
    xmin=0.03125, xmax=1,
    ymin=0.0001, ymax=6,
    xtick={0.03125,0.0625,0.125,0.25,0.5,1},
    legend pos=south east,
    ymajorgrids=true,
    grid style=dashed,
]

\addplot[
    color=blue,
    mark=*,
    ]
    coordinates {
    (1,0.9946) (0.5,0.2870) (0.25,0.0792) (0.125,0.0140) (0.0625,0.0020) (0.03125,0.00054691)
    };

\addplot[
    color=red,
    mark=x,
    ]
    coordinates {
    (1,0.2783) (0.5,0.0561) (0.25,0.0106) (0.125,0.0015) (0.0625,0.0003011) (0.03125,0.00023905)
    };

\addplot[
    color=green,
    mark=square*,
    ]
    coordinates {
    (1,4.5164) (0.5,1.5580) (0.25,0.3994) (0.125,0.1191) (0.0625,0.0301) (0.03125,0.0092)
    };

\addplot[
    color=cyan,
    mark=triangle*,
    ]
    coordinates {
    (1,5.6003) (0.5,2.4957) (0.25,1.1011) (0.125,0.3597) (0.0625,0.0846) (0.03125,0.0192)
    };
    \legend{$L^\infty$, $L^2$, $H^1$, $DG$}
    
 \logLogSlopeTriangle{0.6}{0.1}{0.25}{2}{black};
 \logLogSlopeTriangle{0.6}{0.1}{0.1}{3}{black};
    
\end{loglogaxis}
\end{tikzpicture}
    \caption{Bidomain model, $p=2$}
    \label{fig:bido22}
\end{subfigure} 

\caption{Computed errors vs mesh size (log-log scale) for $p=1,2$: monodomain (top) and bidomain (bottom) models.}
\label{fig:convergence}
\end{figure}
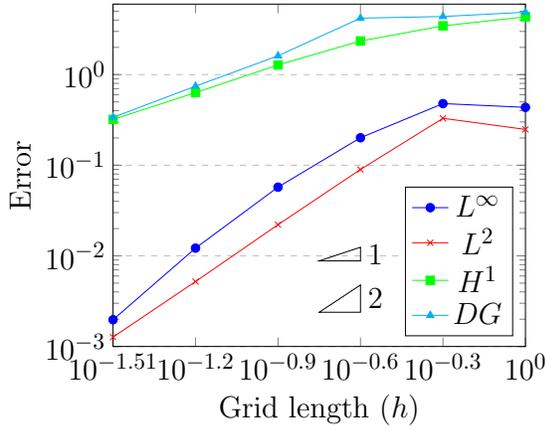
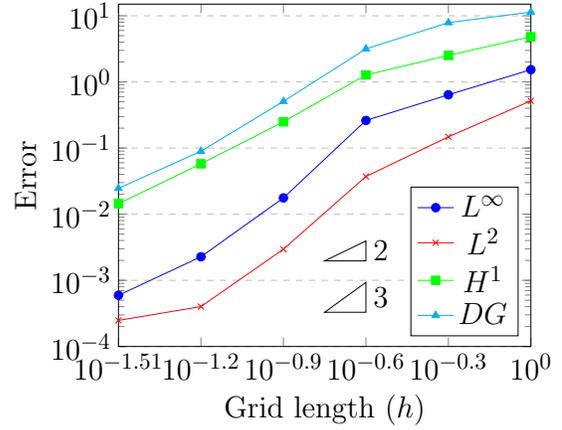
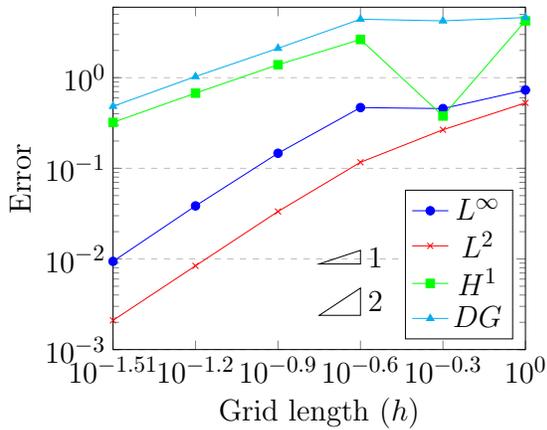
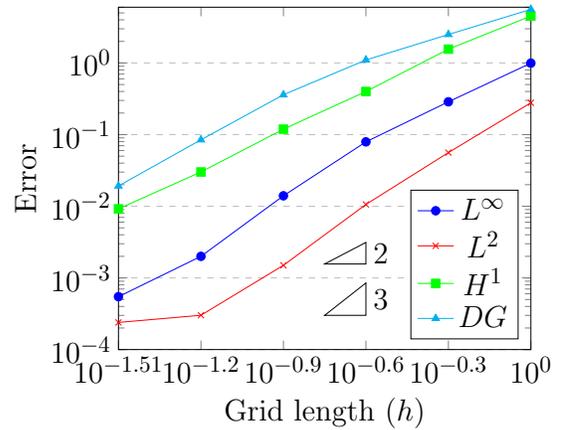

\subsection{Towards a realistic simulation}
 \label{sec:towardsarealisticsimulation}

In this section a first pseudo-realistic simulation is performed in order to assess the coherence between the monodomain model, the numerical discretization and the physical phenomenon. We aim at analysing the behaviour of a localized applied current on a bounded portion of heart tissue. The domain $\Omega \subset \mathbb{R}^2$ is defined as a reference square as in \Cref{sec:verificationtestcases} and boundary conditions are homogeneous so that the only source of potential is represented by $I^{ext}$ while the domain $\Omega$ is isolated from the outer regions. The applied current is defined as 
\begin{equation*}
    I^{ext}(x,y,t)=2\cdot 10^6 \mathcal{I}_{[0.4,0.6]}(x) \mathcal{I}_{[0.4,0.6]}(y) \mathcal{I}_{[0,10^{-3}]}(t),
\end{equation*}
where $\mathcal{I}_{[a,b]}(x)$ is the indicator function, i.e., $\mathcal{I}_{[a,b]}(x)=1$ if $x \in [a,b]$ and 0 otherwise. This definition of $I^{ext}$ represents a temporary and localized electric shock in the center of the domain. Conductivity and ionic model parameters are chosen as in \Cref{sec:verificationtestcases} except for $C_m=10^{-2}$ and $\epsilon=40$. This latter choice allows to obtain electric activation and repolarization in reasonable times. Then, the time step is $\Delta t=10^{-3}$, second order polynomials ($p=2$) are employed as well as a uniformly refined grid with $h=2^{-6}$. 

\begin{figure}
    \centering
    \begin{subfigure}{0.4\textwidth}
            \includegraphics[width=\textwidth]{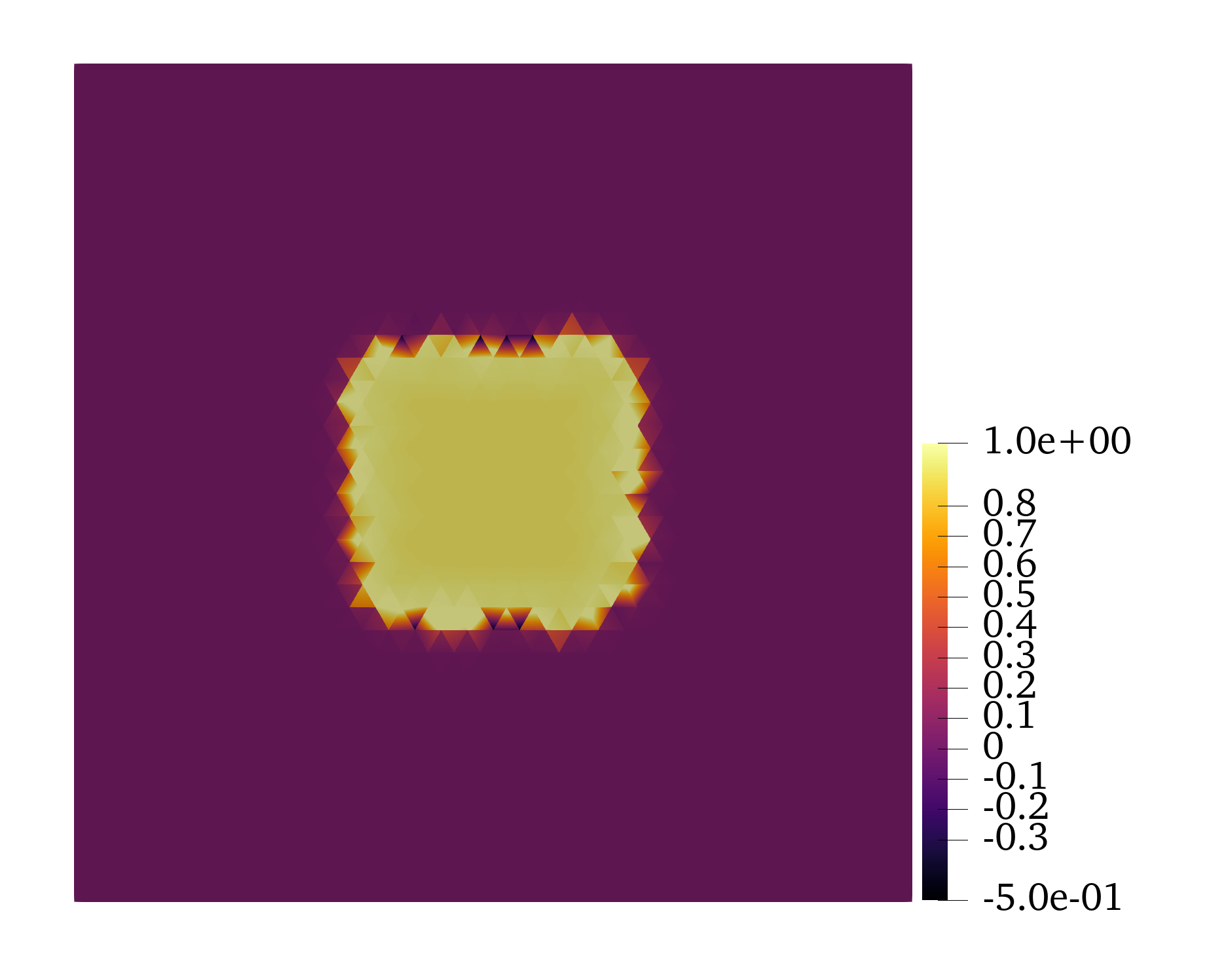}
            \caption{$t=0.04$}
    \end{subfigure}
    \begin{subfigure}{0.4\textwidth}
            \includegraphics[width=\textwidth]{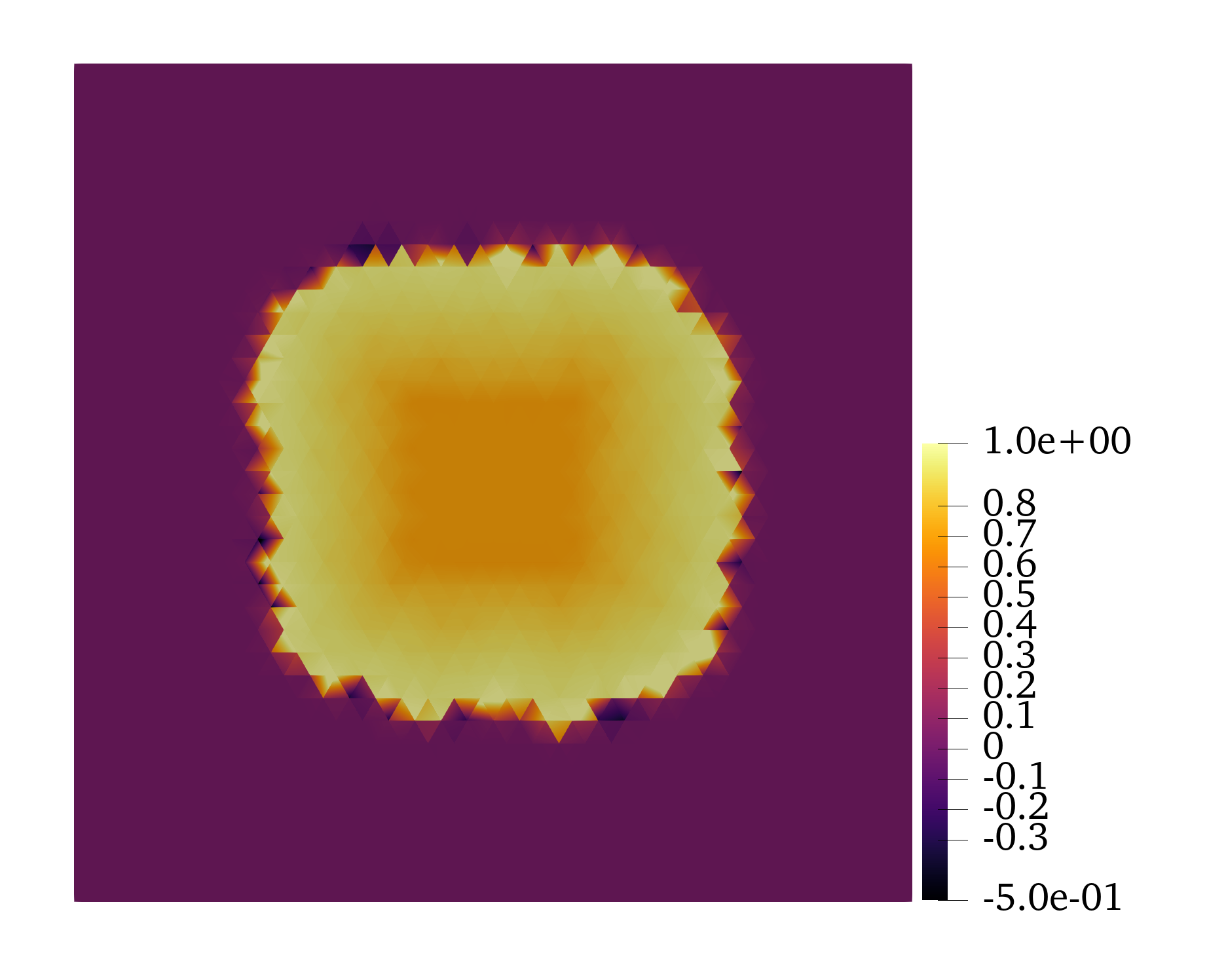}
            \caption{$t=0.10$}
    \end{subfigure}
    \\
    \begin{subfigure}{0.4\textwidth}
            \includegraphics[width=\textwidth]{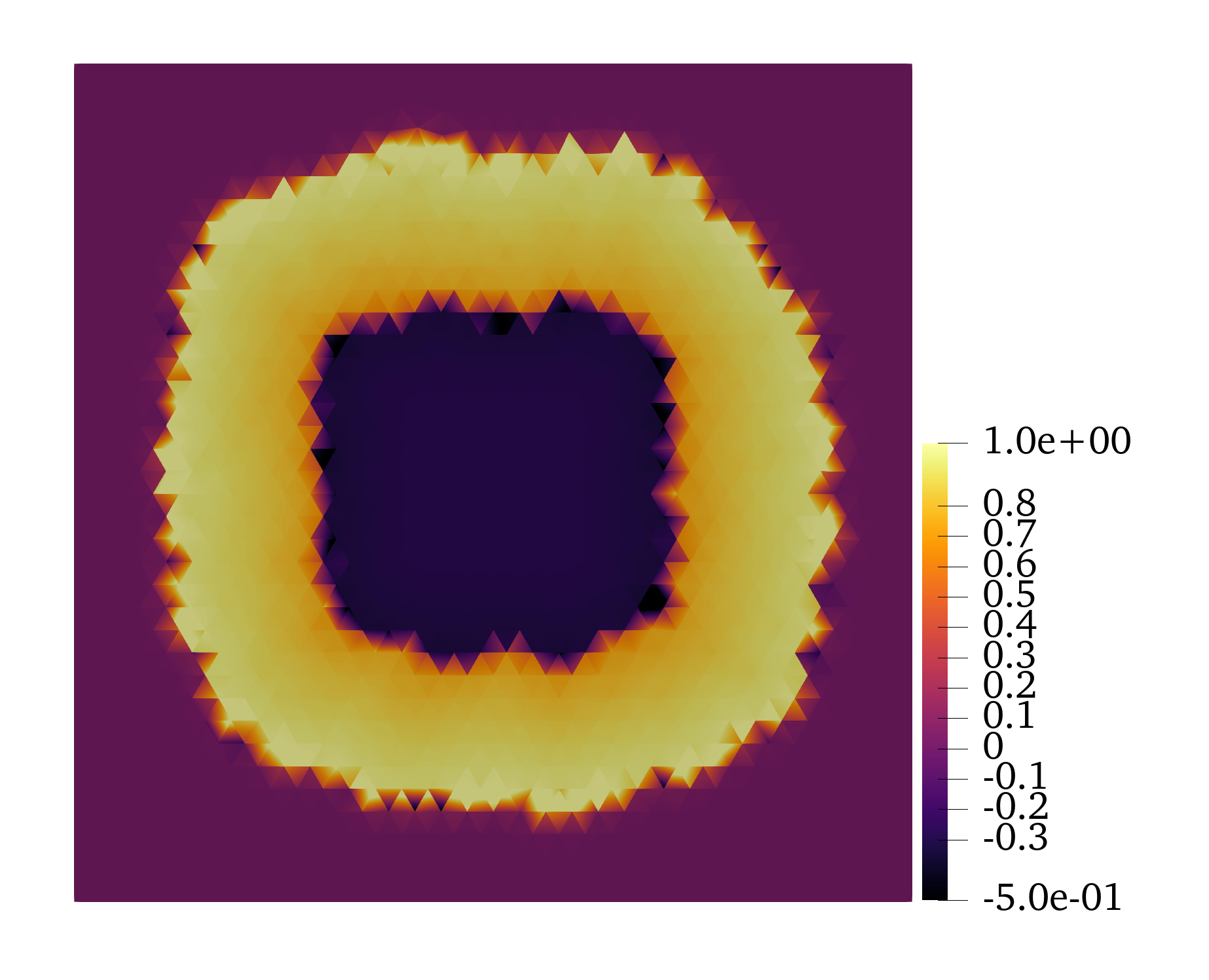}
            \caption{$t=0.16$}
    \end{subfigure}
    \begin{subfigure}{0.4\textwidth}
            \includegraphics[width=\textwidth]{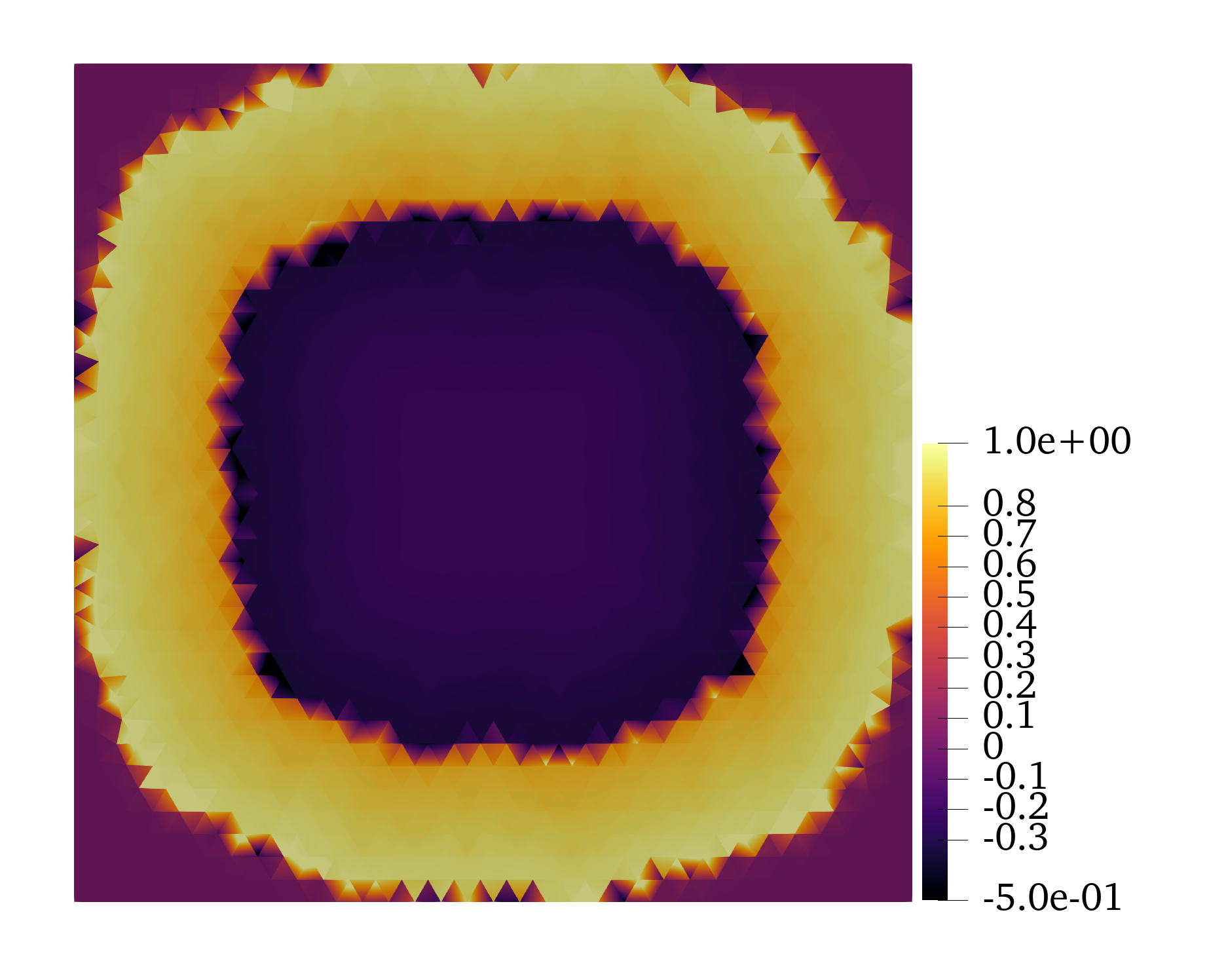}
            \caption{$t=0.22$}
    \end{subfigure}
    \\
    \begin{subfigure}{0.4\textwidth}
            \includegraphics[width=\textwidth]{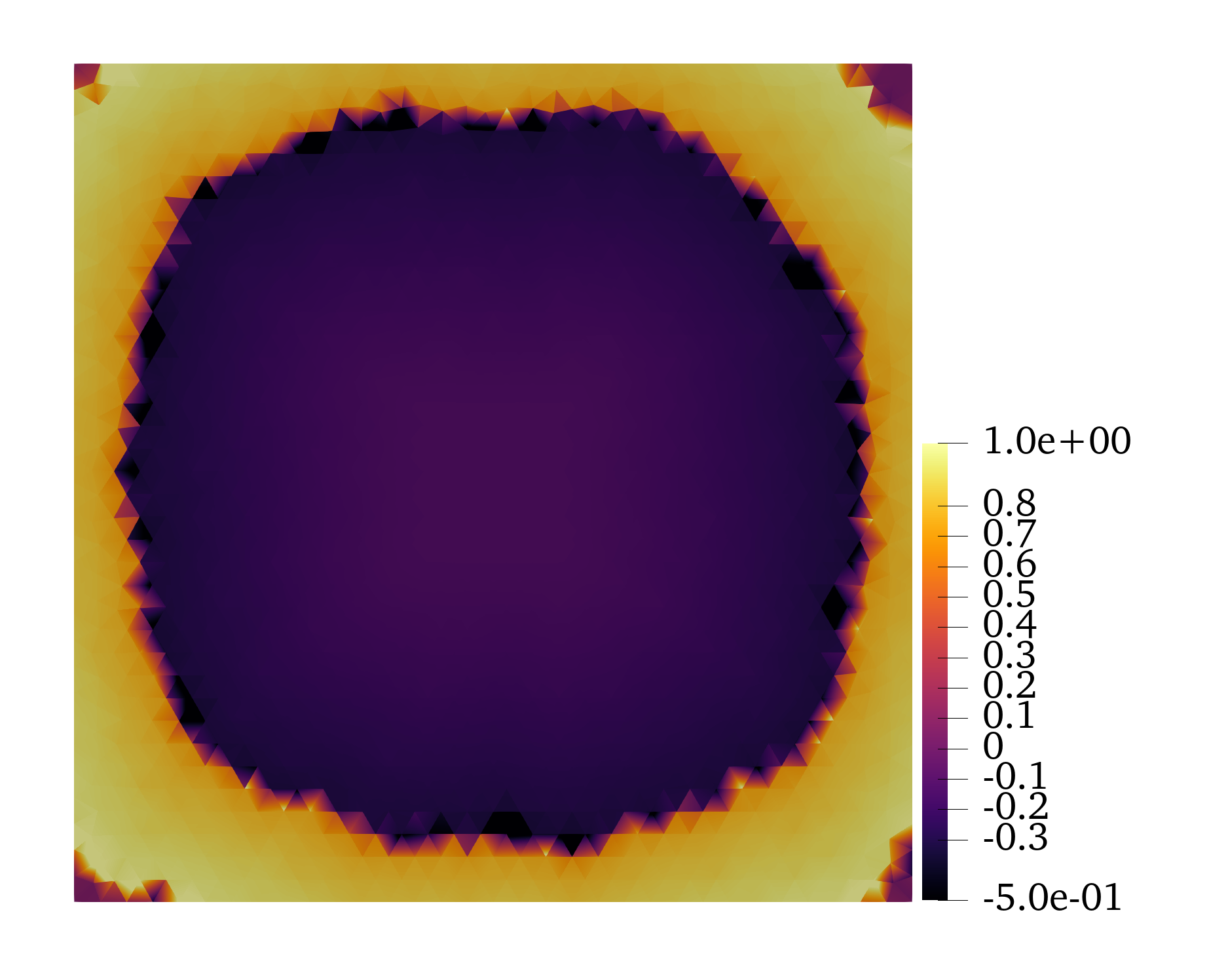}
            \caption{$t=0.28$}
    \end{subfigure}
    \begin{subfigure}{0.4\textwidth}
            \includegraphics[width=\textwidth]{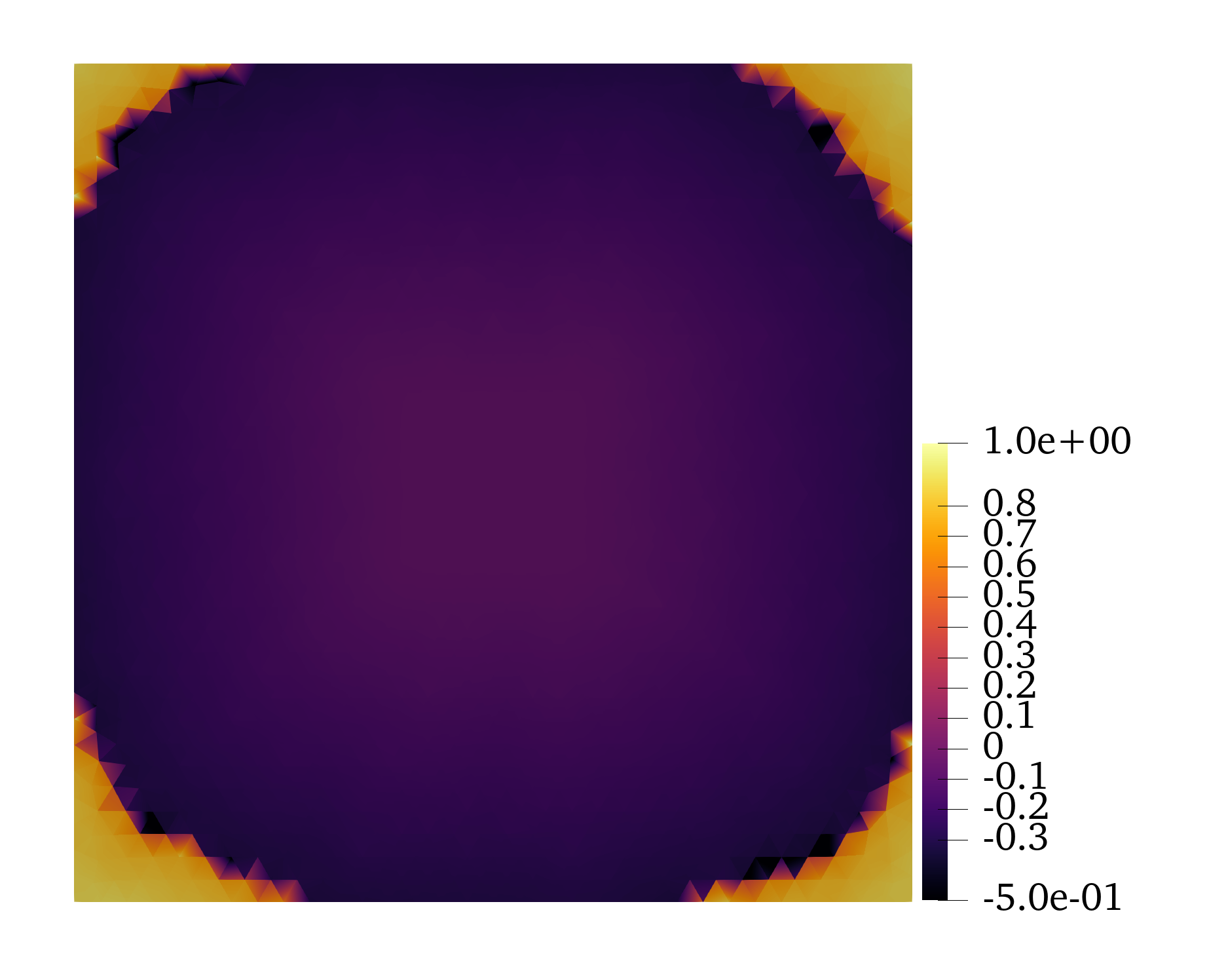}
            \caption{$t=0.34$}
    \end{subfigure}

    \caption{Snapshots from the pseudo-realistic simulation. An external current is applied in the center of the domain initially at rest.}
    \label{fig:realsimulation}
\end{figure}

Finally, results are shown in \Cref{fig:realsimulation} where the activation and subsequent propagation of the electric potential are well visible. The maximum value $V_m=1$ is achieved on the crest while the internal part undergoes repolarization where the potential falls below the rest value and then slowly returns to zero. 

\section{Conclusions}
The numerical discretizations of the monodomain and bidomain models have been successfully carried out through a high-order discontinuous Galerkin formulation and the use of the Dubiner spectral basis. Furthermore, the space discretization has been complemented with a semi-implicit time integration. Finally, the numerical tests in \Cref{sec:numericalresults} have from one side verified the convergence properties of the numerical solution and, on the other hand, compared the numerical simulation with the physiological propagation of the electric potential. In particular, the space discretization error decreases very fast as the grid is refined, meeting the expected theoretical orders of convergence, and soon becoming comparable with the time-step error. Also the physiological test, despite the anisotropy due to the mesh orientation, exhibits good reconstruction of both the depolarization and repolarization phases so that the electric wave propagation is correctly computed.

\section*{Acknowledgements}
\noindent Funding: PCA acknowledges the consortium iNEST (Interconnected North-East Innovation Ecosystem), Piano Nazionale di Ripresa e Resilienza (PNRR) – Missione 4 Componente 2, Investimento 1.5 – D.D. 1058 23/06/2022, ECS00000043, supported by the European Union's NextGenerationEU program and the INdAM - GNCS Project ``Algoritmi efficienti per la gestione e adattazione di mesh poligonali'', codice CUP\_E53C22001930001. PFA has been partially funded by the research grant PRIN2020 n. 20204LN5N5 funded by MUR. PFA has been partially supported by ICSC—Centro Nazionale di Ricerca in High Performance Computing, BigData, and Quantum Computing funded by European Union—NextGeneration EU.
 PCA, PFA, CV are members of the INdAM group GNCS “Gruppo Nazionale per il Calcolo Scientifico” (National Group for Scientific Computing). CV has been partially supported by the Italian Ministry of University and Research (MIUR) within the PRIN (Research projects of
relevant national interest) MIUR PRIN22-PNRR n. P20223KSS2 "Machine learning for fluid-structure interaction in cardiovascular problems: efficient solutions, model reduction, inverse problems, and by  the Italian Ministry of Health within the PNC PROGETTO HUB LIFE SCIENCE - DIAGNOSTICA AVANZATA (HLS-DA) "INNOVA", PNC-E3-2022-23683266–CUP: D43C22004930001, within the "Piano Nazionale Complementare Ecosistema Innovativo della Salute” - Codice univoco investimento: PNC-E3-2022-23683266

\section*{Use of AI tools declaration}
The authors declare they have not used Artificial Intelligence (AI) tools in the creation of this article.

\section*{Conflict of interest}
The authors  declare no conflict of interest.

\bibliographystyle{abbrv}
\bibliography{bibliography.bib}
\end{document}